\documentclass[12pt]{article}
\usepackage{mathrsfs}
\usepackage{amsfonts}
\usepackage{indentfirst}
\usepackage{amssymb}
\usepackage{enumerate}
\usepackage{amsmath}
\usepackage{graphicx}
\usepackage{hyperref}
\usepackage{color}

= 0pt \leftmargin     = 1.25in \topmargin =0pt \headheight     = 0pt \headsep
= 0pt \topskip =0pt
\def\sqr#1#2{{\vcenter{\vbox{\hrule height.#2pt
              \hbox{\vrule width.#2pt height#1pt \kern#1pt \vrule width.#2pt}
              \hrule height.#2pt}}}}
\def\signed #1{{\unskip\nobreak\hfil\penalty50
              \hskip2em\hbox{}\nobreak\hfil#1
              \parfillskip=0pt \finalhyphendemerits=0 \par}}
\def\endpf{\signed {$\sqr69$}}

\def\3n{\negthinspace \negthinspace \negthinspace }
\def\2n{\negthinspace \negthinspace }
\def\1n{\negthinspace }

\def\={\buildrel \triangle \over =}

\def\no{\noindent}

\def\ms{\medskip}
\def\bs{\bigskip}

\def\|{\Big |}
\def\({\Big (}
\def\){\Big )}
\def\[{\Big[}
\def\]{\Big]}
\def\be{\begin{equation}}
\def\bel{\begin{equation}\label}
\def\ee{\end{equation}}
\def\bea{\begin{eqnarray}}
\def\eea{\end{eqnarray}}
\def\bt{\begin{theorem}}
\def\et{\end{theorem}}
\def\bc{\begin{corollary}}
\def\ec{\end{corollary}}
\def\bl{\begin{lemma}}
\def\el{\end{lemma}}
\def\bp{\begin{proposition}}
\def\ep{\end{proposition}}
\def\br{\begin{remark}}
\def\er{\end{remark}}
\def\ba{\begin{array}}
\def\ea{\end{array}}
\def\bd{\begin{definition}}
\def\ed{\end{definition}}
\newtheorem{lemma}{Lemma}[section]
\newtheorem{remark}{Remark}[section]

\newtheorem{theorem}{Theorem}[section]
\newtheorem{corollary}{Corollary}[section]

\newtheorem{definition}{Definition}[section]

\newtheorem{proposition}{Proposition}[section]

\makeatletter
   
   \@addtoreset{equation}{section}
\makeatother
\begin{document}
	\title{\bf  A class of fuzzy numbers induced by probability density functions and their arithmetic operations}
	\author{ Han Wang\thanks{ School of Mathematical Sciences,
			Beijing Normal University,
			100875 Beijing,
			China.   Email: 18843017677@163.com.\ms
		}  $\;$$\;$
		and
		Chuang Zheng\thanks{ School of Mathematical Sciences,
			Beijing Normal University,
			100875 Beijing,
			China.   Email: chuang.zheng@bnu.edu.cn.\ms (The corresponding author)
	}}
	%\date{}
	\maketitle 
	
	\begin{abstract}
		In this paper we are interested in a class of fuzzy numbers which is uniquely identified by their membership functions. The function space, denoted by $X_{h, p}$, will be constructed by combining a class of nonlinear mappings $h$ (subjective perception) and a class of probability density functions  (PDF) $p$ (objective entity), respectively. Under our assumptions, we prove that there always exists a class of $h$ to fulfill the observed outcome for a given class of $p$.  Especially, we prove that the common triangular number can be interpreted by a function pair $(h, p)$.
		As an example, we consider a sample function space $X_{h, p}$ where $h$ is the tangent function and $p$ is chosen as the Gaussian kernel with free variable $\mu$. 
		By means of the free variable $\mu$ (which is also the expectation of $p(x; \mu)$),  we define the addition, scalar multiplication and subtraction on  $X_{h, p}$. We claim that, under our definitions, $X_{h, p}$ has a linear algebra.  Some  numerical examples are provided to illustrate the proposed approach.
	\end{abstract}
	
	\bs

	\no{\bf Key Words}. Fuzzy numbers; Basic concepts; Probability density function;  Gaussian kernel; Fuzzy arithmetic;  Gaussian probability density membership function (G-PDMF)
	
	\section{Introduction}
	
	Fuzzy numbers and fuzzy set theory are topics originated from Zadeh (\cite{Zadeh1965}) by dealing with the imprecise quantities and uncertainty. Since then, they have been  successfully applied in a wide area of topics from pure and applied mathematics, computer science and other related fields, such as fuzzy logic, fuzzy information, soft computing, fuzzy control, etc. 
	
	In general,  a fuzzy number can be uniquely determined by its membership function. In many applications, the membership functions of fuzzy numbers are based on subjective perceptions rather than data or other objective entities involved. The construction of  an appropriate membership function is the cornerstone upon which fuzzy set theory has evolved. See, for instance, \cite{dombi1990, krusinska1986note, 1996Course, wu2013decomposition, wu2020arithmetic, mashchenko2021sums} and the references therein. 
	
	In this paper, we consider a class of fuzzy numbers, denoted by $X_{h, p}$, in which the membership function is carried out by the transformation of the probability density function $p$ combining a nonlinear mapping $h$. Some basic assumptions will be made on the pair $(h,p)$ such as continuity and monotonicity. The exact definition of  the space $X_{h, p}$ will be stated in \eqref{PDMFS}.  We prove that, under the limited information of the fuzzy number, there exists at least one pair $(h,p)$ such that the corresponding membership function fulfills the given data. The details of the description are put in Theorem \ref{1+1} and \ref{m+n}, respectively. Note that comparing to some of the existing methods for obtaining the membership function we listed in the preliminary, one of the advantages of our methodology  is that it includes the subjective factor as well as the objective information. More precisely,  we determine the membership function in two steps: we first put the subjective perception on the type of the pair $(h,p)$ with undetermined parameters. Secondly, by means of the information in the given data, we objectively determine the parameters and identify the exclusive membership function of the fuzzy number. 
	
	Once the membership functions of the fuzzy numbers are identified, one of the basic issues is how to perform the arithmetic operations on them.   In the fuzzy world, the arithmetic operations on real numbers in the classical crisp set level turns to be the algebraic operations on the membership functions.  
	There have been amount of papers studying the fundamental definitions of the operations and corresponding algebraical structures. 
	For instance, triangular norms is introduced to concern the binary operations on the  interval $[0,1]$(see, for instance, \cite{KLEMENT20045}).  
    The interactive  fuzzy numbers and corresponding arithmetic operators are proposed by means of joint possibility distributions (\cite{Zadeh1975199, Dubois1981TAC, Esmi2019246}). Especially, the proposed addition and subtraction can reach the minimum norm compared to other mothod given by sup-J extension principle (\cite{ESMI2021}).  
    The intuitionistic fuzzy set and the corresponding probabilistic addition is constructed and has been successfully handled the fuzzy aggregation problem for expert systems  (\cite{ATANASSOV198687,atanassov1999intuitionistic, SZMIDT2000505, xu2007}).  
	Furthermore, when one considers the fuzzy differential equation,  a common definition of (generalized-) Hukuhara differentiability is required, in which the difference between two fuzzy numbers has to be designed by means of $\alpha$-cut of the corresponding  membership functions (\cite{dubois1982towards, hukuhara1967integration, puri1983differentials}). 
	
	To specify the class of functions under consideration, we fix $h$ as the tangent function and $p$ is given by the Gaussian kernel $p(\cdot, \mu)$, i.e. with $\sigma=1$ and $\mu$ to be undetermined (see the exact definition in Formula  \eqref{abcmumu}). Note that they both fulfill those assumptions in Definitions \ref{h-}--\ref{pdf} as we stated in the preliminaries.  
	We call the function as the form in \eqref{abcmumu} the Gaussian Probability Density Membership Function (G-PDMF) and the corresponding functional space as  G-PDMF Space.  
	
	We design the arithmetic operations, such as addition, scalar multiplication and subtraction on the class of G-PDMFs by means of the expectation parameter $\mu$. The arithmetic operations between fuzzy numbers are transforming to the arithmetic operations between the corresponding parameters $\mu$. As we shall see in Definition \ref{def} and Theorem \ref{linearAlgebra},  the advantage of our design is that we introduce a linear structure on the G-PDMF Space via $\mu$ and $\alpha$-cuts representation is not needed during the computational implemetation. Our work also can be seen as an attempt to make a bridge between probability and fuzzy theory.
	
	The paper is organized as follows. In section \ref{sec2}, the basic concept on fuzzy number and some requirements of the membership function are given. In section  \ref{sec3}, we establish a class of membership function space by introducing a nonlinear map $h$ and the probability density function $p$,  along with the fact that they fulfill the demands of fuzzy numbers.  We also give a constructive proof to show that the common triangular number can be described by our methodology.  In section  \ref{sec4}, we define a sample function space $X_{h,p}$ and address the arithmetic operations on the  G-PDMF. In section  \ref{sec5}, some numerical examples and corresponding graphs are shown to illustrate the operations under consideration. In section \ref{sec6}, we present a final remark to make a complete summary of the paper.

\section{Preliminaries}\label{sec2}
The formulation of membership functions is the crucial step in the design of fuzzy system. There are several methods to develop them. We summarize some of them as follows:
\begin{enumerate}[1)]
	\item L-R linear functions, which is the simplest possible model (\cite{dubois1980fuzzy}); 
	\item Rational functions of polynomials (\cite{giachetti1997parametric, guerra2005approximate});
	\item B-Spline MF (\cite{wang1995fuzzy}); 
	\item Piecewise linear functions (\cite{wen2019novel} and refs [8-20] in it).
\end{enumerate}
In all of these definitions of fuzzy numbers, the membership function needs to satisfy the following assumptions:
\begin{definition}\label{1st_Def}% (\cite{SHEN2020})
	A fuzzy number $A$ is a fuzzy subset of the real line $\mathbb{R}$ with membership function $f_{A}$ which possesses  the following properties: 
	\begin{enumerate}[a)]
		\item $f_{A}$ is fuzzy convex,
		\item $f_{A}$ is normal i.e., $ \exists x_{0} \in \mathbb{R} $ such that $f_{A}(x_{0})=1$,
		\item $f_{A}$ is upper semi-continuous,
		\item The closure of the set $\{x\in\mathbb{R}|f_{A}(x) > 0\}$ is compact.
	\end{enumerate} 
\end{definition}

Definition \ref{1st_Def} is  straightforward and has been used extensively in practical applications  (\cite{SHEN2020}). However, the above conditions are too vague  and a particular class of functions, named as monotonic fuzzy numbers, is introduced  with the  following more precise assumptions:

\begin{definition}\label{2nd_Def}
	A monotonic fuzzy number $\tilde{b}$, denoted by $\tilde{b}=(a,b,c)$,  is defined as a membership function $f(x)$  which possesses  the following properties (\cite{dubois1980fuzzy} ): 
	\begin{enumerate}[a)]
		\item $f(x)$ is increasing on the interval $[a,b]$ and decreasing on $[b,c]$,
		\item $f(x)=1$ for $x=b$, $f(x)=0$ for $x \leq a$ or $x \geq c $,
		\item $f (x)$ is upper semi-continuous,
	\end{enumerate} 
where $a,b,c,$ are real numbers satisfying $-\infty < a\leq b\leq c < +\infty$. 
\end{definition}
Clearly, a class of triangular fuzzy numbers is a subset of the class of monotonic fuzzy numbers. It is due to the fact that, in the definition of the triangular fuzzy numbers, the function in the condition $a)$ of Definition \ref{2nd_Def} is restricted by linear ones  (\cite{liang2013}).

Now we construct a function space containing  membership functions satisfying all requirements above. 

To start with, we first define a nonlinear mapping $h^{-}$ from $[a,b]$ to $\mathbb{R}$, which is crucial to describe the fuzzy number $\tilde{b}=(a, b, c)$. 
\begin{definition}\label{h-}
	Let $h^{-}$ be a function defined on the interval $ (a,b)$. We say $h^{-}$ is a \textbf{left auxiliary function  (LAF)} of the fuzzy number $\tilde{b}$,  if $h^{-}$ satisfies
		\begin{enumerate}[a)]
		\item $\lim_{x \rightarrow a^+}h^{-}(x)=-\infty, \quad \lim_{x \rightarrow b^-}h^{-}(x)=+\infty$,
		\item $h^{-}$ is continuous on $ (a,b)$,
		\item $h^{-}$ is increasing on $ (a,b)$.
	\end{enumerate} 
\end{definition}
Similarly, we define $h^{+}$ on the right side $[b,c]$ as follows:
\begin{definition}\label{h+}
We say  $h^{+}$ is a \textbf{right auxiliary function (RAF) } of $\tilde{b}$,  if $h^{+}: (b,c)\rightarrow\mathbb{R}$ satisfies
\begin{enumerate}[a)]
	\item $\lim_{x \rightarrow b^+}h^{+}(x)=+\infty, \quad \lim_{x \rightarrow c^-}h^{+}(x)=-\infty$,
	\item $h^{+}$ is continuous on $(b,c)$,
	\item $h^{+}$ is decreasing on $(b,c)$.
\end{enumerate} 
\begin{figure}[htbp]
	\centering
	\begin{minipage}[t]{0.48\textwidth}
		\centering
		\includegraphics[width=6cm]{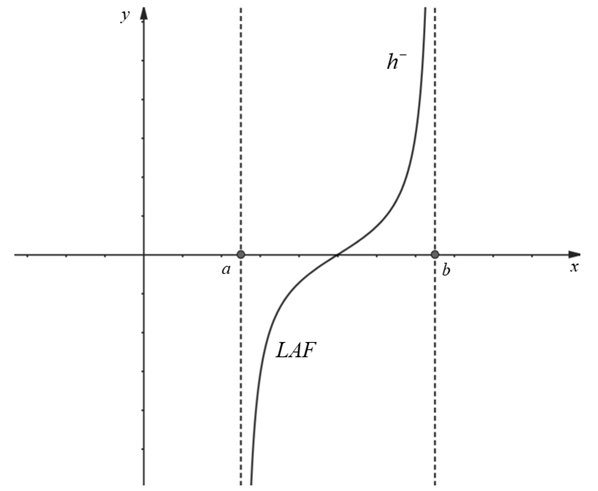}
		\caption{Shape of the LAF}\label{Fig.LAF} 
	\end{minipage}
	\begin{minipage}[t]{0.48\textwidth}
		\centering
		\includegraphics[width=6cm]{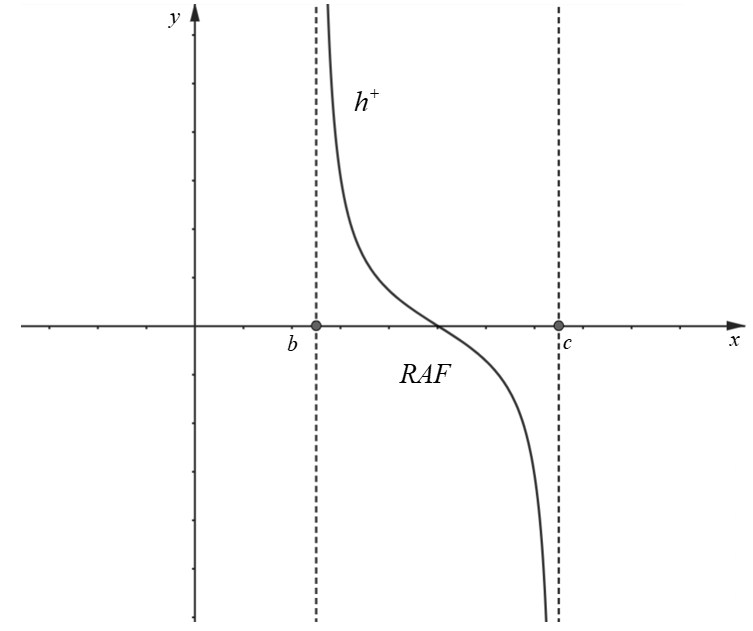}
		\caption{Shape of the RAF}\label{Fig.RAF} 
	\end{minipage}
\end{figure}
\end{definition}

The probability density function is defined by
\begin{definition}\label{pdf}
	We say  $p^-$ and $p^+$ are \textbf{probability density functions (PDFs)},  if $p^-$ and $p^+$ both  satisfy
	\begin{enumerate}[a)]
		\item $\displaystyle\int_{-\infty}^{+\infty}p^-(y)dy=1,\;\int_{-\infty}^{+\infty}p^+(y)dy=1,$
		\item$p^-(t) \geq 0, \; p^+(t) \geq 0,\qquad \forall t \in (-\infty,+\infty)$.
	\end{enumerate} 
\end{definition}

Note that to construct the desired memebership function, we need two different PDFs on intervals $(a, b)$ and $(b,c)$. In fact, $p^-$ is the PDF used on $(a, b)$ and $p^+$ is the PDF used on $(b,c)$, respectively.

Based on above  functions $h^{-}, h^{+}$ and $p^{-},p^{+}$, the fuzzy number $\tilde{b}$ is constructed by
\begin{definition}\label{PMF}
We say $\tilde{b}$ is a  fuzzy number generated by $h^{-}, h^{+}$ and $p^{-},p^{+}$,  if the membership function of $\tilde{b}$ has the form
\begin{equation}\label{fxabc}
\begin{array}{lll}	
f(x)%&=&\chi_{[a,b]}(x)f_{-}(x)+\chi_{(b,c]}(x)f_{+}(x)	\\	
&=&\left \{
\begin{aligned}
&0, &x \in (-\infty,a]\\
&\int_{-\infty}^{h^{-}(x)}p^{-}(y)dy,  &x\in (a,b)\\
& 1 , &x=b\\
&\int_{-\infty}^{h^{+}(x)}p^{+}(y)dy,  &x\in (b,c)\\
&0,    &x\in [c,+\infty)\\
\end{aligned}
\right.
\end{array}
\end{equation} 
\end{definition}
Note that $h^{-}$ and $h^{+}$ can be derived from the same function $h\in H([0,1])$ with
\begin{equation}\label{spaceH}
H([0,1])=\{
h \| h \; \hbox{ is a LAF on } [0,1]
\},
\end{equation}
i.e. $h$ satisfies Definition \ref{h-} with $a=0$ and $b=1$. For simplicity, in the sequel, we only consider the case that $h^{-}$ and $h^{+}$ is originated from the same $h\in H ([0,1])$, i.e., 
$$
h^{-}(x)=h(\frac{x-a}{b-a}), \qquad 
h^{+}(x)=h(\frac{c-x}{c-b}).
$$
Similarly, we take $p^-$ and $p^+$ is originated from the same class of PDFs.

Let $h\in H([0,1])$ and $p=(p^-,p^+)$ satisfy Definition \ref{pdf}, we call the function space 
\begin{equation}\label{PDMFS}
	X_{h,p}:=\{f_{h,p}(x):\mathbb{R} \rightarrow [0,1]  \;\; \hbox{is as the form of \eqref{fxabc}} :a\leq b\leq c\}
\end{equation}
a \textbf{Probability Density Membership Function} (PDMF) space.

To sum up, we have:
\begin{theorem}\label{PDMF} 
The PDMF in the space  $X_{h,p}$  as the form in \eqref{fxabc}  fulfills all requirements in Definition \ref{1st_Def} and Definition \ref{2nd_Def}.
\end{theorem}
{\bf Proof:} It is straightforward.

\section{PDMFs with control points}\label{sec3}

In this section we consider the case that some control points are predetermined on the shape of the membership function $f_{h,p}\in X_{h,p}$.
\begin{theorem}\label{1+1}
	Let 
	\begin{equation}\label{PQ}
		P=(x^-, y^-)\in (a,b)\times(0,1),\qquad Q=(x^+, y^+)\in (b,c)\times(0,1).
	\end{equation} 
Then  there exists  at least one pair $(h, p)$ such that the graph of $f_{h,p}$ passes through $P$ and $Q$, i.e., $f_{h,p}(x^-)=y^-$ and $f_{h,p}(x^+)=y^+$. 
\end{theorem}
{\bf Proof:}
We first fix a function $h\in H([0,1])$ satisfying  \eqref{spaceH}. Set 
\begin{equation*}
	z^{-}=h^{-}(x^{-})=h(\dfrac{x^{-}-a}{b-a}), \qquad z^{+}=h^{+}(x^{+})=h(\dfrac{c-x^{+}}{c-b})
\end{equation*}  
Thus $-\infty<z^{-},z^{+}<+\infty$. 
We now construct functions $p^{-},p^{+}$ from $\mathbb{R}$ to $\mathbb{R}^+\cup\{0\}$ as follows: 
\begin{equation}\nonumber
	p^{-}(t)=\left\{
	\begin{aligned}
		&0,&t\in(-\infty,z^{-}-1]\\
		&y^{-},&t\in(z^{-}-1,z^{-}]\\
		&1-y^{-},&t\in(z^{-},z^{-}+1]\\
		&0,&t\in(z^{-}+1,+\infty)\\	
	\end{aligned}
	\right.
\end{equation}

\begin{equation}\nonumber
	p^{+}(t)=\left\{
	\begin{aligned}
		&0,&t\in(-\infty,z^{+}-1]\\
		&y^{+},&t\in(z^{+}-1,z^{+}]\\
		&1-y^{+},&t\in(z^{+},z^{+}+1]\\
		&0,&t\in(z^{+}+1,+\infty)\\	
	\end{aligned}
	\right.
\end{equation}
By direct computation, we have $ p^{-}(t),p^{+}(t)\geq 0 ,\forall t \in \mathbb{R} $ and
\begin{equation*}
	\int_{-\infty}^{+\infty}p^{-}(t)dt=1,	\int_{-\infty}^{z^{-}}p^{-}(t)dt=y^{-},\\
	\int_{-\infty}^{+\infty}p^{+}(t)dt=1,	\int_{-\infty}^{z^{+}}p^{+}(t)dt=y^{+}.\\
\end{equation*}
Hence  $p^{-},p^{+}$  satisfy Definition \ref{pdf} and 
\begin{equation*}
	\int_{-\infty}^{h^{-}(x^{-})}p^{-}(t)dt=y^{-},
	\int_{-\infty}^{h^{+}(x^{+})}p^{+}(t)dt=y^{+},
\end{equation*}
which means  $f_{h,p}(x^{-})=y^{-}$ and $f_{h,p}(x^{+})=y^{+}$  with  $h=(h^{-},h^{+}),p=(p^{-},p^{+})$.
 This completes the proof.

Consequently, we have the following theorem for finite number of control points.
\begin{theorem}\label{m+n}
 Let $P_i=(x_i^-, y_i^-)\in (a,b)\times(0,1), i=1,\cdots, m$ satisfy $a<x_1^-\cdots<x_m^-<b$ and $0<y_1^-\cdots<y_m^-<1$, $Q_j=(x_j^+, y_j^+)\in (b,c)\times(0,1),  j=1,\cdots, n$ satisfy $b<x_1^+\cdots<x_n^-<c$ and $0<y_{n}^{+}<\cdots<y_{1}^{+}<1$. Then there exists  at least one pair $(h, p)$ such that the graph of $f_{h,p}$ passes through all points $P_i$'s and $Q_j$'s, i.e., $f_{h,p}(x_i^-)=y_i^-, i=1,\cdots, m$ and $f_{h,p}(x_j^+)=y_j^+, j=1,\cdots, n$. 
\end{theorem}
{\bf Proof:}
We first fix a function $h\in H([0,1])$ satisfying  \eqref{spaceH}. Set   
\begin{center}
	$z_{i}^{-}=h(\frac{x_{i}^{-}-a}{b-a}), \qquad i=1,\cdots,m. $
\end{center}
The monotony of $h$ indicates that  $z_{i}^{-}<z_{i'}^{-}$ for all $ 1\leq i<i'\leq m$.
We now construct a function $p^{-}:\mathbb{R}\rightarrow \mathbb{R}^+\cup\{0\}$ as follows:  
\begin{equation*} 
	p^{-}(t)=\left\{
	\begin{aligned}
		&0,&t\in(-\infty,z_{1}^{-}-1]\\
		&y_{1}^{-},&t\in(z_{1}^{-}-1,z_{1}^{-}]\\
		&\dfrac{y_{2}^{-}-y_{1}^{-}}{z_{2}^{-}-z_{1}^{-}},&t\in(z_{1}^{-},z_{2}^{-}]\\
		&\dfrac{y_{3}^{-}-y_{2}^{-}}{z_{3}^{-}-z_{2}^{-}},&t\in(z_{2}^{-},z_{3}^{-}]\\
		&\cdots,\\
		&\dfrac{y_{m}^{-}-y_{m-1}^{-}}{z_{m}^{-}-z_{m-1}^{-}},&t\in(z_{m-1}^{-},z_{m}^{-}]\\
		&1-y_{m}^{-},&t \in (z_{m}^{-},z_{m}^{-}+1]\\
		&0,&t\in(z_{m}^{-}+1,+\infty)
	\end{aligned}
	\right.
\end{equation*}
By direct computation, we have $ p^{-}(t)\geq 0,\forall t \in \mathbb{R} $ and
\begin{equation*}
	\int_{-\infty}^{+\infty}p^{-}(t)dt=1,
	\int_{-\infty}^{z_{i}^{-}}p^{-}(t)dt=y_{i}^{-}.
\end{equation*}
So the function $p^{-}$ satisfies Definition \ref{pdf} and 
\begin{equation*}
	\int_{-\infty}^{h(\frac{x_{i}^{-}-a}{b-a})}p^{-}(t)dt=y^{-}_{i},\qquad i=1,\cdots,m.
\end{equation*}
The similar result can be proven for points $Q_{j}(j=1,\cdots,n)$ and the function $h$. In fact, by the same procedure, we can construct a function $p^{+}$ satisfying Definition \ref{pdf} and 
\begin{equation*}
	\int_{-\infty}^{h(\frac{c-x^{+}_{j}}{c-b})}p^{+}(t)dt=y^{+}_{j},\qquad j=1,\cdots,n.
\end{equation*}
Hence,  $f_{h,p}(x_i^{-})=y_i^{-} (i=1,\cdots,m)$ and $f_{h,p}(x_j^{+})=y_j^{+} (j=1,\cdots,n)$.
This completes the proof.\endpf

Note that the membership function we construct there can be seen as a similar  form of pentagon fuzzy numbers with $m=n=1$ (\cite{mondal2017pentagonal}) and B-spline fuzzy numbers with $n+m$ control points (\cite{wang1995fuzzy}), respectively.

As a direct consequence of our result,  we have following theorem:

\begin{theorem} There exists 
	at least one pair $(h,p)$ such that the triangular fuzzy number $(a,b,c)$ is in the PDMF space $X_{h,p}$.\\
\end{theorem}	

	\textbf{Proof:} 
	Recall that a triangular fuzzy number determined by the triplet $(a,b,c)$ of real numbers with $a<b<c$ has a membership function as follows:
	\begin{equation}\nonumber
		(a,b,c)=\left\{
		\begin{aligned}
			&0,&x\in(-\infty,a]\\
			&\dfrac{x-a}{b-a},&x\in(a,b)\\
			&1,&x=b\\
			&\dfrac{c-x}{c-b},&x\in(b,c)\\
			&0,&x\in[c,+\infty)
		\end{aligned}
		\right.
	\end{equation}
	
	\noindent We first fix $\mu \in \mathbb{R}$ and  $p^{-}(t)=p^{+}(t)=p(t; \mu)=\frac{1}{\sqrt{ 2\pi}}e^{-\frac12(t-\mu)^{2}}$.
	 Notice that the standard normal cumulative distribution function 
		$F(x)=\int_{-\infty}^{x}\frac{1}{\sqrt{2\pi}}e^{-\frac{t^{2}}{2}}dt$  is a strictly monotonic function on $\mathbb{R}$. According to the inverse function theorem, the inverse function of $F$ exists, i.e. its quantile function $Q(y)=\inf\{x\in\mathbb{R}:y\leq F(x)\}$ exists. 
 Hence, for all $x^{-}\in (a,b)$, we can set $z^{-}$ satisfying the equation
	\begin{equation*}
		\int_{-\infty}^{z^{-}}p(t;\mu)dt=\int_{-\infty}^{z^{-}-\mu^{-}}\frac{1}{\sqrt{2\pi}}e^{-\frac{t^{2}}{2}}dt=\frac{x^{-}-a}{b-a}.
	\end{equation*}
    Since $0<\frac{x^{-}-a}{b-a} < 1$, we can define a function $h:(0,1)\rightarrow (-\infty,+\infty)$ such that 
	\begin{equation*}
		h(\frac{x^{-}-a}{b-a})=z^{-}.\\
	\end{equation*} 
	By direct computation, we have 
	\begin{equation*}
		\int_{-\infty}^{h(\frac{x^{-}-a}{b-a})}p(t)dt=\frac{x^{-}-a}{b-a}.
	\end{equation*}
	and 
	for all $ x^{+} \in (b,c),$
	\begin{equation*}
		\int_{-\infty}^{h(\frac{c-x^{+}}{c-b})}p(t)dt=\frac{c-x^{+}}{c-b}.
	\end{equation*}

We  now verify that $h$ is a LAF as in Definition \ref{h-}.  In fact, we have
	\begin{enumerate}[\bf{Claim} a)]
		\item   $h(0)=-\infty,h(1)=+\infty$. \\
		 \textbf{Proof:}Set $h(0)=z_{0}^{-}$, i.e.
		\begin{equation*}
			\int_{-\infty}^{z_{0}^{-}}\frac{1}{\sqrt{ 2\pi}}e^{-\frac12(t-\mu)^{2}}dt=0.
		\end{equation*}
		Consequently $z_{0}^{-}=-\infty$. Similarly we can prove $h(1)=+\infty$.
		\item  $h$ is increasing.\\
		 \textbf{Proof:}
		For all $ x_{1}^{-},x_{2}^{-}\in (a,b),x_{1}^{-} < x_{2}^{-}$, we suppose that $h(x_{i}^{-})=z_{i}^{-}(i=1,2)$ which means
		\begin{equation*}
			\int_{-\infty}^{z_{1}^{-}}\frac{1}{\sqrt{2\pi}}e^{-\frac{(t-\mu)^{2}}{2}}dt=\frac{x_{1}^{-}-a}{b-a},\quad
			\int_{-\infty}^{z_{2}^{-}}\frac{1}{\sqrt{2\pi}}e^{-\frac{(t-\mu)^{2}}{2}}dt=\frac{x_{2}^{-}-a}{b-a}.
		\end{equation*}
		After linear transformation, we have
		\begin{equation*}
			\int_{-\infty}^{z_{1}^{-}-\mu}\frac{1}{\sqrt{2\pi}}e^{-\frac{t^{2}}{2}}dt=\frac{x_{1}^{-}-a}{b-a},\quad
			\int_{-\infty}^{z_{2}^{-}-\mu}\frac{1}{\sqrt{2\pi}}e^{-\frac{t^{2}}{2}}dt=\frac{x_{2}^{-}-a}{b-a}.
		\end{equation*}
		The monotony of standard normal cumulative distribution function indicates that $z_{1}^{-}<z_{2}^{-}$, i.e. $h(x_{1}^{-})<h(x_{2}^{-})$.
		\item $h$ is a continuous function.\\
		 \textbf{Proof:}
		We first fix $ x_{0}^{-} \in (a,b). $ Then according to inverse function theorem, there exists a $z_{0}^{-}$ which satisfies $z_{0}^{-}=h(x_{0}^{-})$. For  all $\varepsilon >0$, we set $z_{0}^{-}-\varepsilon=h(x_{1}^{-})$ and $z_{0}^{-}+\varepsilon=h(x_{2}^{-})$, then $\exists$ $\delta=\min\{x_{0}^{-}-x_{1}^{-},x_{2}^{-}-x_{0}^{-}\} >0,$ s.t. if $|x-x_{0}^{-}|<\delta$,  we have $|z-z_{0}^{-}|<\varepsilon$.
		
	\end{enumerate}

	Consequently the function $h$ we constructed above satisfies Definition \ref{h-}
    and $f_{h,p}(x)=(a,b,c)$. This completes the proof. \endpf
	
\section{Gaussian PDMFs with $\sigma=1$}\label{sec4}

In this section we  first establish a sample space  by taking $h$ as a tangent function and $p$ as a Gaussian membership function with the standard deviation $\sigma=1$, respectively. Two control points $P$ and $Q$ will be given on the shape of the membership function $f_{h,p}(x)$. The expectation $\mu$ of the Gaussian Kernel will be determined by the control points $P$ and $Q$. Consequently, by means of the parameter $\mu$, we design the operations on $f_{h,p}(x)$ such as addition, scalar multiplication and subtraction. Some properties and advantages of our definitions are given. 
\subsection{Definitions}\label{sec4.1}
Set
\begin{equation}\label{hp}
h(x) =\tan(\pi x - \frac{\pi}{2}),  \; x \in (0,1), 
\qquad p(t)=p (t; \mu)=\frac{1}{\sqrt{ 2\pi}}e^{-\frac12(t-\mu)^{2}}, \; t\in\mathbb{R}.
\end{equation} 
Denoted by 
\begin{equation*}
\varphi(x;a,b) =\dfrac{x-a}{b-a}, \; x \in [a,b]; 
\qquad
\varphi(x;b,c) =\dfrac{x-b}{c-b},   \; x \in [b,c], 
\end{equation*} 
the LMF and RMF are given by
\begin{equation}\label{LRMF}
\begin{aligned}
h^{-}(x)&= h(\varphi(x;a,b))=\tan(\frac{\pi}{b-a}(x-a)-\frac{\pi}{2}), \quad &x\in (a, b),   \\
h^{+}(x)&=h(1-\varphi(x;b,c))=\tan(\frac{\pi}{c-b}(c-x)-\frac{\pi}{2}), \quad &x\in(b, c) .
\end{aligned}
\end{equation}
Moreover, we assume that there are two control points $P(x^{-}, y^{-}), Q(x^{+}, y^{+})$ on each side of the central value $b$ with $a<x^{-}<b<x^{+}<c$. The corresponding membership function as in \eqref{fxabc} has the exact form
\begin{equation}\label{abcmumu}
f_{h,p}(x)=
\left \{
\begin{aligned}	
&0, &x \in (-\infty,a]\\
&f_{-}(x;\mu^{-},a,b),  &x\in (a,b)\\
& 1 , &x=b\\
&f_{+}(x;\mu^{+},b,c),  &x\in (b,c)\\
&0,    &x\in [c,+\infty)\\
\end{aligned}
\right.
\end{equation}
where $f_-$ and $f_+$ is given by
\begin{equation}\nonumber
\begin{aligned}
	f_{-}(x;\mu^{-},a,b)
	&=\int_{-\infty}^{h^{-}(x)}p(t;\mu^{-})dt    & \\ 
	&=\int_{-\infty}^{\tan(\frac{\pi}{b-a}(x-a)-\frac{\pi}{2})}\frac{1}{\sqrt{ 2\pi}}e^{-\frac12(t-\mu^{-})^{2}}dt, &\qquad x\in(a,b)\\ 
\end{aligned}
\end{equation}
\begin{equation}\nonumber
	\begin{aligned}
		f_{+}(x;\mu^{+},b,c)
		&=\int_{-\infty}^{h^{+}(x)}p(t;\mu^{+})dt   & \\ 
		&=\int_{-\infty}^{\tan(\frac{\pi}{c-b}(c-x)-\frac{\pi}{2})}\frac{1}{\sqrt{ 2\pi}}e^{-\frac12(t-\mu^{+})^{2}}dt, &\qquad x\in(b,c).\\ 
	\end{aligned}
\end{equation}

 The following Theorem holds:
\begin{theorem}\label{abcPQ}
	The function space $X_{h,p}$ is a PDMFS as the form in \eqref{PDMFS} if $p, h$ are taken as  in \eqref{hp}. Moreover, there exists a unique pair $(\mu^{-},\mu^{+})$ for any $(P, Q)$ given by \eqref{PQ}.
	We call the above PDMFS as a Gaussian PDMFS, abbreviated by G-PDMFS. 
\end{theorem}

Some remarks are in order.

\begin{remark}\label{h-p}
	Note that the pair $(h,p)$ we designed as in \eqref{hp} is the subjective perception we  offer to the class of the membership functions under consideration. The parameter $\mu^-$ (resp. $\mu^+$) is  remained to be uniquely determined by the pre-given information $P(x^-,y^-)$ (resp. $Q(x^+, y^+)$). We emphasize that, rather than the tangent function, there are plenty of possibilities to choose $h$ as a pre-designed function,  such as  Logit function  or inverse sigmoid function, etc.
\end{remark}
\begin{remark}\label{G-PDMFS}
	As a direct consequence of Theorem \ref{abcPQ}, it is reasonable to give two equivalent notations of the G-PDMF as  
	\begin{equation}\label{abcde}
	 \langle(a,b,c);P,Q\rangle
	\Longleftrightarrow
	\langle(a,b,c);\mu^{-},\mu^{+}\rangle. 
	\end{equation}
	As we shall see in the proof,  $(\mu^{-},\mu^{+})$ can be uniquely identified by means of the inverse function of Formula (\ref{mumu}).
\end{remark}
\begin{remark}
	Mathematically speaking,  ranking fuzzy numbers can be seen as ranking the functions as the form of \eqref{abcmumu} in G-PDMF space.  Various definitions of ranking methods can be designed based on either of the two equivalent notations in \eqref{abcde}. For instance, $b_1\preceq b_2$ if $ (\mu^-_1)^2+(\mu^+_1)^2\leq(\mu^-_2)^2+(\mu^+_2)^2$, with $b_i=\langle(a,b,c);\mu^{-}_i,\mu^{+}_i\rangle, i=1,2$. The detail of the design is beyond the scope of this paper and will  be discussed elsewhere.
\end{remark}

{\bf Proof:}
For  $(h,p)$ given in \eqref{hp},  it is obvious that  $h$ belongs to $H([0,1])$ and  $p$ is a PDF. We only need to verify that $(\mu^-, \mu^+)$ is uniquely determined by $(P, Q)$. In fact, set $z^{-}=h^{-}(x^{-})=h(\dfrac{x-a}{b-a})$, it follows that
\begin{equation}\label{mumu}
	\int_{-\infty}^{z^{-}}p(t;\mu^{-})dt=\int_{-\infty}^{z^{-}}\frac{1}{\sqrt{2\pi}}e^{-\frac{(t-\mu^{-})^{2}}{2}}dt=y^{-}.
\end{equation}
According to inverse function theorem, the inverse function of standard normal cumulative distribution function exists. Thus there must exist a $\mu^{-}$ satisfying the equation above. To verify the $\mu^{-}$ is unique, we suppose that there exists two values $\mu_{1}$ and $\mu_{2}$ satisfying
\begin{equation*}
	\int_{-\infty}^{z^{-}}\frac{1}{\sqrt{2\pi}}e^{-\frac{(t-\mu_{1})^{2}}{2}}dt=y^{-}, \quad
	\int_{-\infty}^{z^{-}}\frac{1}{\sqrt{2\pi}}e^{-\frac{(t-\mu_{2})^{2}}{2}}dt=y^{-}.
\end{equation*}  
Consequently, 
\begin{equation*}
	\int_{-\infty}^{z^{-}-\mu_{1}}\frac{1}{\sqrt{2\pi}}e^{-\frac{t^{2}}{2}}dt=y^{-},\quad
	\int_{-\infty}^{z^{-}-\mu_{2}}\frac{1}{\sqrt{2\pi}}e^{-\frac{t^{2}}{2}}dt=y^{-}.
\end{equation*}
The monotony of standard normal cumulative distribution function indicates that $\mu_{1}=\mu_{2}$. The proof for $\mu^{+}$ is similar and we omit it.

\subsection{Operational laws}
For $h$ and the PDF $p$  given by \eqref{hp}, we now design operational laws of the G-PDMFS
\begin{equation}\label{GPDMFS}
X_{h, p}(\mathbb{R})
=
\{\tilde b\|\tilde b=\langle(a,b,c);\mu^{-},\mu^{+}\rangle \;\hbox{has the form of \eqref{abcmumu}} \},
\end{equation}
 such as addition, scalar multiplication and subtraction. 

\begin{definition} \label{def}Let $\tilde{b}_{1}=\langle(a_{1},b_{1},c_{1});\mu^{-}_{1},\mu^{+}_{1} \rangle$ and $\tilde{b}_{2}=\langle(a_{2},b_{2},c_{2});\mu^{-}_{2},\mu^{+}_{2} \rangle$ be two G-PDMFs in $X_{h,p}$, then
	\begin{enumerate}[(1)]
	\item  $\tilde{b}_{1}\oplus\tilde{b}_{2}=\langle(a_{1}+a_{2},b_{1}+b_{2},c_{1}+c_{2});\mu^{-}_{1}+\mu^{-}_{2},\mu^{+}_{1}+\mu^{+}_{2} \rangle.$
	\item 
	$\lambda \tilde{b}_{1}=\left\{
	\begin{array}{l}
	\langle(\lambda a_{1},\lambda b_{1},\lambda c_{1});\lambda \mu^{-}_{1},\lambda\mu^{+}_{1} \rangle, \forall \lambda \geq 0.\\
	\langle(\lambda c_{1},\lambda b_{1},\lambda a_{1});\lambda \mu^{+}_{1},\lambda\mu^{-}_{1} \rangle, \forall \lambda <  0.\\
	\end{array}
	\right.$
	\item
	$\tilde{b}_{1}\ominus\tilde{b}_{2}=\langle(a_{1}-c_{2},b_{1}-b_{2},c_{1}-a_{2});\mu^{-}_{1}-\mu^{+}_{2},\mu^{+}_{1}-\mu^{-}_{2} \rangle.$
\end{enumerate}
\end{definition}

Some remarks are in order:
\begin{remark}
	To adapt the rules of the standard arithmetic addition on real numbers,  it is mandatory to require that,  for $\tilde{b}_{3}=\tilde{b}_{1}\oplus\tilde{b}_{2}$, the membership function of $\tilde{b}_{3}$ satisfies $f(b_1+b_2)=1$. Nevertheless, the definitions on the endpoints $a_3, c_3$ can be varied  in several ways. In fact, it is reasonable to choose any  $a_3, b_3$ such that 
	\begin{itemize}
		\item $a_{3} \in [a_{1}+a_{2},b_{1}+b_{2}-min\{b_{1}-a_{1},b_{2}-a_{2}\}] $,
		\item  $c_{3} \in [b_{1}+b_{2}+min\{c_{1}-b_{1},c_{2}-b_{2}\}, c_{1}+c_{2}]$,
	\end{itemize}
	depending on the real-world situations for the fuzzy system.  We speculate that different choices  of endpoints may leads to diverse properties.  The same methodology holds true on the design of scalar multiplication and subtraction. Especially, to obey the rule that $a\leq b\leq c$ for any G-PDMF $\langle(a,b,c);\mu^{-},\mu^{+}\rangle$, it is  natural to design  the scalar multiplication separately for different sign of $\lambda$ as in $(2)$ of Definition \ref{def}.
\end{remark}
\begin{remark}
	Note that it is straightforward to define the subtraction on $X_{h, p}$ as a generalization of its scalar multiplication with $\lambda =-1$ (see Formula ($3$) of  Definition \ref{def}). Hence, the class of G-PMDFs has a linear algebra generated by the expectation $\mu$ of the Gaussian Kernel. We believe that this kind of definitions between two G-PMDFs can be further exploited to some kind of fuzzy differentiation and integration. Consequently, one can construct corresponding fuzzy differential equations and further works need to be done.
\end{remark}
\begin{remark}\label{implementation}
    Note that according to Formula ($3$), the subtraction of two G-PDMFs is straightforward since $(\mu^{-},\mu^{+})$ can be uniquely identified via the inverse function of Formula (\ref{mumu}), as we claimed in Remark \ref{G-PDMFS}. More precisely, set two G-PDMFs $\tilde{b}_{1}$ and $\tilde{b}_{2}$ with the information 
    \begin{equation}\label{b1b24}
    \tilde{b}_{1}=\langle(a_1,b_1,c_1);P_1,Q_1\rangle,\;
    \tilde{b}_{2}=\langle(a_2,b_2,c_2);P_2,Q_2\rangle.
    \end{equation}
    Their subtraction $\tilde{b}_{1}\ominus\tilde{b}_{2}$, denoted by $\tilde{X}$,  will be computed by two steps:
    \begin{enumerate}
        \item By means of the inverse function of Formula (\ref{mumu}), we compute the corresponding expectation values $\mu$ and \eqref{b1b24} turns to 
        $$
         \tilde{b}_{1}=\langle(a_1,b_1,c_1);\mu_1^-,\mu_1^+\rangle,\;
        \tilde{b}_{2}=\langle(a_2,b_2,c_2);\mu_2^-,\mu_2^+\rangle.
        $$
        \item  We compute $\tilde{X}$ directly by Formula ($3$).
    \end{enumerate}
    With this procedure, one can easily construct the method of solving the toy fuzzy equation $\tilde{b}_{2}\oplus \tilde{X}=\tilde{b}_{1}$ with direct computation $\tilde{X}=\tilde{b}_{1}\ominus\tilde{b}_{2}$. A specific example is also stated in Section \ref{sec5}.
\end{remark}
\begin{remark}
	Note that the operations via $\alpha$-cut is not needed in our design with the numerical implementation in Remark \ref{implementation}.  Frankly speaking,  the initial datum is given by $\langle(a, b,c); P, Q\rangle$ with $P=(x^-, y^-)$ and $Q=(x^+, y^+)$, which is the partial information of the fuzzy number. Notice that we do not know the value of $f(x)$ at other points other that $x=a, b, c, x^-, x^+$ on the interval $[a,c]$.  Once we obtain the parameters $(\mu^-, \mu^+)$, the shape of the membership function is fixed and the corresponding arithmetic operations between fuzzy numbers are transforming to the arithmetic operations of $(\mu^-, \mu^+)$ and we have all information of the fuzzy number.
\end{remark}

The following properties hold:
\begin{theorem} Let $\tilde{b}_{1}=\langle(a_{1},b_{1},c_{1});\mu^{-}_{1},\mu^{+}_{1}\rangle$, and $\tilde{b}_{2}=\langle(a_{2},b_{2},c_{2});\mu^{-}_{2},\mu^{+}_{2}\rangle$ be two G-PDMFs in $X_{h,p}$, then there exist  G-PDMFs $\tilde{b}_{3},\tilde{b}_{4},\tilde{b}_{5} \in X_{h,p}$  such that  $ \tilde{b}_{3}=\tilde{b}_{1}\oplus\tilde{b}_{2}$ , $\tilde{b}_{4}=\tilde{b}_{1}\ominus\tilde{b}_{2}$ and $\tilde{b}_{5}=\lambda \tilde{b}_{1}(\forall \lambda \in \mathbb{R})$, respectively.
\end{theorem}
{\bf Proof:}
Since $\tilde{b}_{1},\tilde{b}_{2}$ are two G-PDMFs, we have $a_{1},b_{1},c_{1},\mu^{-}_{1},\mu^{+}_{1},a_{2},b_{2},c_{2},\mu^{-}_{2},\mu^{+}_{2}\in \mathbb{R}$. Thus $a_{1}+a_{2},b_{1}+b_{2},c_{1}+c_{2},a_{1}-c_{2},b_{1}-b_{2},c_{1}-a_{2},\mu^{-}_{1}+\mu^{-}_{2},\mu^{+}_{1}+\mu^{+}_{2},\mu^{-}_{1}-\mu^{+}_{2},\mu^{+}_{1}-\mu^{-}_{2},\lambda a_{1},\lambda b_{1},\lambda c_{1},\lambda \mu^{-}_{1},\lambda \mu^{+}_{1} \in \mathbb{R}$.
It follows from Definition \ref{def} that 
\begin{equation}\nonumber
	\begin{aligned}
		\tilde{b}_{3}&=\langle(a_{1}+a_{2},b_{1}+b_{2},c_{1}+c_{2});\mu^{-}_{1}+\mu^{-}_{2},\mu^{+}_{1}+\mu^{+}_{2}\rangle,	\\
		\tilde{b}_{4}&=\langle(a_{1}-c_{2},b_{1}-b_{2},c_{1}-a_{2});\mu^{-}_{1}-\mu^{+}_{2},\mu^{+}_{1}-\mu^{-}_{2} \rangle,\\
		\tilde{b}_{5}&=\left\{
		\begin{array}{l}
			\langle(\lambda a_{1},\lambda b_{1},\lambda c_{1});\lambda \mu^{-}_{1},\lambda\mu^{+}_{1} \rangle, \forall \lambda \geq 0.\\
			\langle(\lambda c_{1},\lambda b_{1},\lambda a_{1});\lambda \mu^{+}_{1},\lambda\mu^{-}_{1} \rangle, \forall \lambda < 0.\\
		\end{array}
		\right.\\ 	
	\end{aligned}
\end{equation}
It is easy to verify that $\tilde{b}_{3},\tilde{b}_{4},\tilde{b}_{5}$ are also G-PDMFs. 
\endpf
\begin{theorem}\label{linearAlgebra}
	Let $\tilde{b}_{i}=\langle(a_{i},b_{i},c_{i});\mu^{-}_{i},\mu^{+}_{i}\rangle, i=1,2,3$ be three G-PDMFs in $X_{h,p}$, then for all $\lambda,\lambda_{1},\lambda_{2} \in \mathbb{R}$, we have 
	\begin{enumerate}[(1)]
		\item $\tilde{b}_{1}\oplus\tilde{b}_{2}=\tilde{b}_{2}\oplus\tilde{b}_{1},(\tilde{b}_{1}\oplus\tilde{b}_{2})\oplus \tilde{b}_{3}=\tilde{b}_{1}\oplus(\tilde{b}_{2}\oplus\tilde{b}_{3})$; 
		\item $\tilde{b}_{1}\ominus\tilde{b}_{2}=-(\tilde{b}_{2}\ominus\tilde{b}_{1})$,  $(\tilde{b}_{1}\ominus\tilde{b}_{2})\ominus \tilde{b}_{3}=\tilde{b}_{1}\ominus(\tilde{b}_{2}\oplus\tilde{b}_{3})$; 
		\item $\lambda(\tilde{b}_{1}\oplus\tilde{b}_{2})=\lambda\tilde{b}_{1}\oplus\lambda \tilde{b}_{2}$, for all $\lambda \in \mathbb{R}$; 
		\item  $\lambda_{1}\tilde{b}_{1} \oplus \lambda_{2}\tilde{b}_{1}= (\lambda_{1}+\lambda_{2})\tilde{b}_{1},$ for all $\lambda_{1},\lambda_{2} \in \mathbb{R}$ and $\lambda_{1}\lambda_{2}\geq 0$;
		\item $\tilde{b}_{1}\oplus(-1)\tilde{b}_{2}=\tilde{b}_{1}\ominus\tilde{b}_{2}.$
	\end{enumerate}
\end{theorem}
{\bf Proof:} Assertion $(1)$  is trivial.

For $(2)$, by the operational law $(3)$ in Definition \ref{def}, we have 
\begin{equation}\nonumber
\begin{aligned}
\tilde{b}_{1}\ominus\tilde{b}_{2}&=\langle(a_{1}-c_{2},b_{1}-b_{2},c_{1}-a_{2});\mu^{-}_{1}-\mu^{+}_{2},\mu^{+}_{1}-\mu^{-}_{2} \rangle,     \\
\tilde{b}_{2}\ominus\tilde{b}_{1}&=\langle(a_{2}-c_{1},b_{2}-b_{1},c_{2}-a_{1});\mu^{-}_{2}-\mu^{+}_{1},\mu^{+}_{2}-\mu^{-}_{1} \rangle.    
\end{aligned}
\end{equation}
Then, by the operational law $(2)$ in Definition \ref{def}, it follows that 
\begin{equation*}
\begin{aligned}
-(\tilde{b}_{2}\ominus\tilde{b}_{1})
&=\langle(-(c_{2}-a_{1}),-(b_{2}-b_{1}),-(a_{2}-c_{1});-(\mu^{+}_{2}-\mu^{-}_{1} ),-(\mu^{-}_{2}-\mu^{+}_{1})\rangle.\\
&=\langle(a_{1}-c_{2},b_{1}-b_{2},c_{1}-a_{2});\mu^{-}_{1}-\mu^{+}_{2},\mu^{+}_{1}-\mu^{-}_{2} \rangle.
\end{aligned}
\end{equation*}
Hence
\begin{equation*}
\tilde{b}_{1}\ominus\tilde{b}_{2}=-(\tilde{b}_{2}\ominus\tilde{b}_{1}).
\end{equation*}	
Also since
\begin{equation*}
\begin{aligned}
(\tilde{b}_{1}\ominus\tilde{b}_{2})\ominus \tilde{b}_{3}&=
\langle(a_{1}-c_{2}-c_{3},b_{1}-b_{2}-b_{3},c_{1}-a_{2}-a_{3});\mu^{-}_{1}-\mu^{+}_{2}-\mu_{3}^{+},\mu^{+}_{1}-\mu^{-}_{2}-\mu_{3}^{-} \rangle	
\end{aligned}		
\end{equation*}	
and
\begin{equation*}
\tilde{b}_{2}\oplus\tilde{b}_{3}=\langle(a_{2}+a_{3},b_{2}+b_{3},c_{2}+c_{3});\mu^{-}_{2}+\mu^{-}_{3},\mu^{+}_{2}+\mu^{+}_{3} \rangle,
\end{equation*}		
then
\begin{equation*}
\begin{aligned}
\tilde{b}_{1}\ominus(\tilde{b}_{2}\oplus\tilde{b}_{3})&=
\langle(a_{1}-(c_{2}+c_{3}),b_{1}-(b_{2}+b_{3}),c_{1}-(a_{2}+a_{3}));\mu_{1}^{-}-(\mu^{+}_{2}+\mu^{+}_{3}),\mu_{1}^{+}-(\mu^{-}_{2}+\mu^{-}_{3})\rangle \\
&=\langle (a_{1}-c_{2}-c_{3},b_{1}-b_{2}-b_{3},c_{1}-a_{2}-a_{3});\mu_{1}^{-}-\mu^{+}_{2}-\mu^{+}_{3},\mu_{1}^{+}-\mu^{-}_{2}-\mu^{-}_{3})\rangle \\
&=	(\tilde{b}_{1}\ominus\tilde{b}_{2})\ominus \tilde{b}_{3}.
\end{aligned}
\end{equation*}

For $(3)$, first we assume that $\lambda \geq 0$.
Then, by the operational law as in Definition \ref{def}, it follows that 
\begin{equation*}
\lambda(\tilde{b}_{1}\oplus\tilde{b}_{2})=\langle(\lambda(a_{1}+a_{2}),\lambda(b_{1}+b_{2}),\lambda(c_{1}+c_{2}));\lambda(\mu^{-}_{1}+\mu^{-}_{2}),\lambda(\mu^{+}_{1}+\mu^{+}_{2})\rangle.
\end{equation*}
Since 
\begin{equation*}
\lambda\tilde{b}_{1}=\langle(\lambda a_{1},\lambda b_{1},\lambda c_{1});\lambda \mu^{-}_{1},\lambda \mu^{+}_{1}\rangle, \lambda\tilde{b}_{2}=\langle(\lambda a_{2},\lambda b_{2},\lambda c_{2});\lambda \mu^{-}_{2},\lambda \mu^{+}_{2}\rangle,
\end{equation*}
then
\begin{equation*}
\lambda\tilde{b}_{1}\oplus\lambda \tilde{b}_{2}=\langle(\lambda a_{1}+\lambda a_{2},\lambda b_{1}+\lambda b_{2},\lambda c_{1}+\lambda c_{2});\lambda \mu^{-}_{1}+\lambda \mu^{-}_{2},\lambda \mu^{+}_{1}+\lambda \mu^{+}_{2}\rangle,
\end{equation*}
hence 
\begin{equation*}
\lambda(\tilde{b}_{1}\oplus\tilde{b}_{2})=\lambda\tilde{b}_{1}\oplus\lambda \tilde{b}_{2}.
\end{equation*}
For the case $\lambda<0$, we have
\begin{equation*}
\begin{aligned}
\lambda(\tilde{b}_{1}\oplus\tilde{b}_{2})&=\langle(\lambda(c_{1}+c_{2})),\lambda(b_{1}+b_{2}),\lambda(a_{1}+a_{2}));\lambda(\mu^{+}_{1}+\mu^{+}_{2}),\lambda(\mu^{-}_{1}+\mu^{-}_{2})\rangle,\\
\lambda\tilde{b}_{1}&=\langle(\lambda c_{1},\lambda b_{1},\lambda a_{1});\lambda \mu^{+}_{1},\lambda \mu^{-}_{1}\rangle,\\
\lambda\tilde{b}_{2}&=\langle(\lambda c_{2},\lambda b_{2},\lambda a_{2});\lambda \mu^{+}_{2},\lambda \mu^{-}_{2}\rangle.
\end{aligned}
\end{equation*}
then
\begin{equation*}
\begin{aligned}
\lambda\tilde{b}_{1}\oplus\lambda\tilde{b}_{2}&=\langle(\lambda c_{1}+\lambda c_{2},\lambda b_{1}+\lambda b_{2},\lambda a_{1}+\lambda a_{2});\lambda \mu^{+}_{1}+\lambda \mu^{+}_{2},\lambda \mu^{-}_{1}+\lambda \mu^{-}_{2}\rangle,\\
&=\langle(\lambda(c_{1}+c_{2})),\lambda(b_{1}+b_{2}),\lambda(a_{1}+a_{2}));\lambda(\mu^{+}_{1}+\mu^{+}_{2}),\lambda(\mu^{-}_{1}+\mu^{-}_{2})\rangle, \\
&=\lambda(\tilde{b}_{1}\oplus\tilde{b}_{2}) \\
\end{aligned}
\end{equation*}

For (4), we first assume the case $\lambda_{1}\geq 0,\lambda_{2}\geq 0.$
Since
\begin{equation*}
\lambda_{1}\tilde{b}_{1}=\langle(\lambda_{1} a_{1},\lambda_{1} b_{1},\lambda_{1} c_{1});\lambda_{1} \mu^{-}_{1},\lambda_{1} \mu^{+}_{1}\rangle,\\
\lambda_{2}\tilde{b}_{1}=\langle(\lambda_{2} a_{1},\lambda_{2} b_{1},\lambda_{2} c_{1});\lambda_{2} \mu^{-}_{1},\lambda_{2} \mu^{+}_{1}\rangle,\\
\end{equation*}
then
\begin{equation}\nonumber
\begin{aligned}
\lambda_{1}\tilde{b}_{1} \oplus \lambda_{2}\tilde{b}_{1}&=\langle(\lambda_{1} a_{1}+\lambda_{2} a_{1},\lambda_{1} b_{1}+\lambda_{2} b_{1},\lambda_{1} c_{1}+\lambda_{2} c_{1});\lambda_{1} \mu^{-}_{1}+\lambda_{2} \mu^{-}_{1},\lambda_{1} \mu^{+}_{1}+\lambda_{2} \mu^{+}_{1}\rangle\\
&=\langle((\lambda_{1}+\lambda_{2}) a_{1},(\lambda_{1}+\lambda_{2}) b_{1},(\lambda_{1}+\lambda_{2}) c_{1});(\lambda_{1}+\lambda_{2}) \mu^{-}_{1},(\lambda_{1}+\lambda_{2}) \mu^{+}_{1}\rangle\\
&=(\lambda_{1}+\lambda_{2})\langle(a_{1},b_{1},c_{1});\mu^{-}_{1},\mu^{+}_{1}\rangle\\
&=(\lambda_{1}+\lambda_{2})\tilde{b}_{1}.
\end{aligned} 
\end{equation}
Secondly, for the case $\lambda_{1}\leq 0,\lambda_{2} \leq 0$,
we have
\begin{equation*}
\lambda_{1}\tilde{b}_{1}=\langle(\lambda_{1} c_{1},\lambda_{1} b_{1},\lambda_{1} a_{1});\lambda_{1} \mu^{+}_{1},\lambda_{1} \mu^{-}_{1}\rangle,\\
\lambda_{2}\tilde{b}_{1}=\langle(\lambda_{2} c_{1},\lambda_{2} b_{1},\lambda_{2} a_{1});\lambda_{2} \mu^{+}_{1},\lambda_{2} \mu^{-}_{1}\rangle.\\
\end{equation*}
Then
\begin{equation}\nonumber
\begin{aligned}
\lambda_{1}\tilde{b}_{1} \oplus \lambda_{2}\tilde{b}_{1}&=\langle(\lambda_{1} c_{1}+\lambda_{2} c_{1},\lambda_{1} b_{1}+\lambda_{2} b_{1},\lambda_{1} a_{1}+\lambda_{2} a_{1});\lambda_{1} \mu^{+}_{1}+\lambda_{2} \mu^{+}_{1},\lambda_{1} \mu^{-}_{1}+\lambda_{2} \mu^{-}_{1}\rangle\\
&=\langle((\lambda_{1}+\lambda_{2}) c_{1},(\lambda_{1}+\lambda_{2}) b_{1},(\lambda_{1}+\lambda_{2}) a_{1});(\lambda_{1}+\lambda_{2}) \mu^{+}_{1},(\lambda_{1}+\lambda_{2}) \mu^{-}_{1}\rangle\\
&=(\lambda_{1}+\lambda_{2})\langle(a_{1},b_{1},c_{1});\mu^{-}_{1},\mu^{+}_{1}\rangle\\
&=(\lambda_{1}+\lambda_{2})\tilde{b}_{1}.
\end{aligned} 
\end{equation}

For (5), since 
\begin{equation*}
(-1)\tilde{b}_{2}=\langle(-c_{2},-b_{2},-a_{2});-\mu^{+}_{2},-\mu^{-}_{2}\rangle,
\end{equation*}
then 
\begin{equation*}
\tilde{b}_{1}\oplus(-1)\tilde{b}_{2}=\langle(a_{1}-c_{2},b_{1}-b_{2},c_{1}-a_{2});\mu^{-}_{1}-\mu^{+}_{2},\mu^{+}_{1}-\mu^{-}_{2}\rangle.
\end{equation*}
Also by the operational law $(3)$ in Definition \ref{def}, we have 
\begin{equation*}
\tilde{b}_{1}\ominus\tilde{b}_{2}=\langle(a_{1}-c_{2},b_{1}-b_{2},c_{1}-a_{2});\mu^{-}_{1}-\mu^{+}_{2},\mu^{+}_{1}-\mu^{-}_{2} \rangle.
\end{equation*}
Thus 
\begin{equation*}
\tilde{b}_{1}\oplus(-1)\tilde{b}_{2}=\tilde{b}_{1}\ominus\tilde{b}_{2},
\end{equation*}
which completes the proof. \endpf

\section{Examples}\label{sec5}

In this section we present some visual graphs to illustrate the shape of G-PDMFs and results of their fuzzy arithmetic operations.
Recall that the original data will be given as the form $\langle(a, b, c);P, Q \rangle$ where 
$$ 
P=(x^-, y^-)\in(a,b)\times(0,1),\qquad Q=(x^+,y^+)\in(b,c)\times(0,1).
$$
To facilitate the operations, we  first compute the expectations $(\mu^{-},\mu^{+})$ of the corresponding Gaussian Kernels. It can be done  via Formula (\ref{mumu}) according to Theorem \ref{abcPQ}. In the sequel, we use the notation $\langle(a, b, c);\mu^{-},\mu^{+}\rangle$ to facilitate the operations. 

\textbf{Example 1.} (Addition)
Let  $\tilde{b}_{1}=\langle(-1,0,1);(-0.5,0.5),(0.5,0.5)\rangle$,  $\tilde{b}_{2}=\langle(-1,1,4);(0,0.5), \allowbreak(2.5,0.5)\rangle$ be two Gaussian PDMFs in $X_{h,p}$ with $(h,p)=(\tan, \mathcal{N}(\mu,1))$ as in Subsection \ref{sec4.1}. According to Theorem \ref{abcPQ},  $(\mu^{-},\mu^{+})$ can be computed via Formula (\ref{mumu}).  Hence, $\tilde{b}_{1}=\langle(-1,0,1);0,0\rangle$ and $\tilde{b}_{2}\langle(-1,1,4);0,0 \rangle$.

By means of  the operation law $(1)$ in Definition \ref{def}, we have
\begin{equation*} 
\tilde{b}_{1}\oplus\tilde{b}_{2}=\langle(-1,0,1);0,0\rangle\oplus \langle(-1,1,4);0,0 \rangle =\langle(-2,1,5);0,0 \rangle.
\end{equation*}
\begin{figure}[htbp]
	\centering
	\begin{minipage}[t]{0.48\textwidth}
		\centering
		\includegraphics[width=6cm]{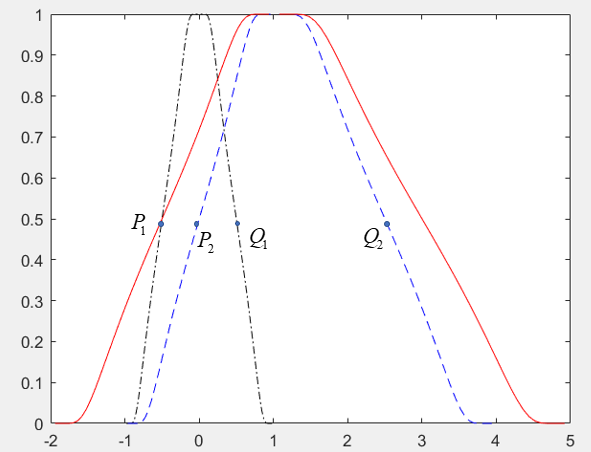}
		\caption{$ \tilde{b}_{1}$ (pecked line), $\tilde{b}_{2}$ (dash line) and $\tilde{b}_{1}\oplus\tilde{b}_{2}$ (solid line)}\label{Fig.1} 
	\end{minipage}
	\begin{minipage}[t]{0.48\textwidth}
		\centering
		\includegraphics[width=6cm]{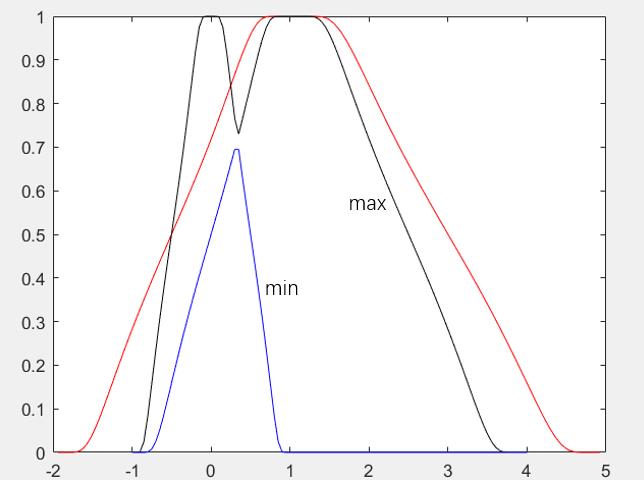}
		\caption{Three results based on MIN, MAX and Formula $(1)$ in Definition \ref{def}}\label{Fig.2} 
	\end{minipage}
\end{figure}

Figure \ref{Fig.1} shows the shape of $ \tilde{b}_{1},\tilde{b}_{2}$ and $\tilde{b}_{1}\oplus\tilde{b}_{2}$, respectively. Figure \ref{Fig.2} presents the fuzzy additions on the two Gaussian PDMFs $\tilde{b}_{1},\tilde{b}_{2}$ using  MAX,  MIN and our operation law $(1)$ in Definition \ref{def}. In Figure \ref{Fig.2},  MIN produces the  lowest membership value for the resulting fuzzy number for all $x\in\mathbb{R}$.  Interestingly, our design produces partially larger value than the one caused by the MAX operation, which may be useful in the realistic application.

\textbf{Example 2.} (Scalar multiplication) 
Let $\tilde{b}_{3}=\langle(-1,1,2);(0,0.75),(1.5,0.6)\rangle$ be a Gaussian PDMF in 
$X_{h,p}$. As above, we can approximately rewrite $\tilde{b}_{3}$ as $\langle(-1,1,2);-0.6745,-0.4399\rangle$.
The Gaussian PDMF resulting from fuzzy scalar multiplication is calculated using the operation law $(2)$ in Definition \ref{def}:
\begin{equation*} 
\begin{aligned}
3\;\tilde{b}_{3}&=3\;\langle(-1,1,2);-0.6745,-0.4399\rangle=\langle(-3,3,6);-2.0225,-1.3197\rangle.\\
(-2)\;\tilde{b}_{3}&=
(-2)\langle(-1,1,2);-0.6745,-0.4399\rangle=\langle(-4,-2,2);0.8798,1.3490\rangle
\end{aligned}\\
\end{equation*}
Figure \ref{Fig.3} and \ref{Fig.4} show the $3 \tilde{b}_{3}$ and $(-2) \tilde{b}_{3}$, respectively:
\begin{figure}[htbp]
	\centering
	\begin{minipage}[t]{0.48\textwidth}
		\centering
		\includegraphics[width=6cm]{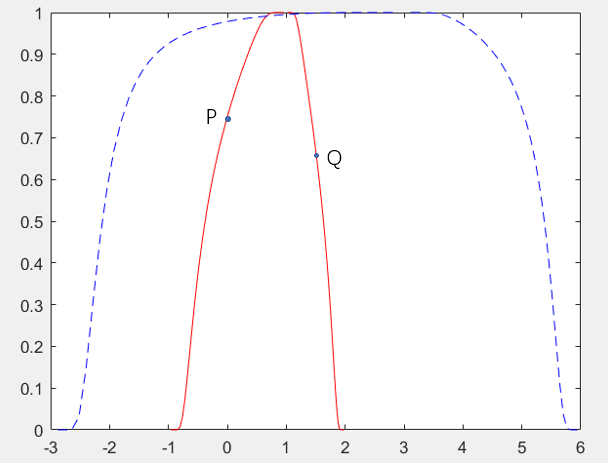}
		\caption{The graph of $3\tilde{b}_{3}$}\label{Fig.3} 
	\end{minipage}
	\begin{minipage}[t]{0.48\textwidth}
		\centering
		\includegraphics[width=6cm]{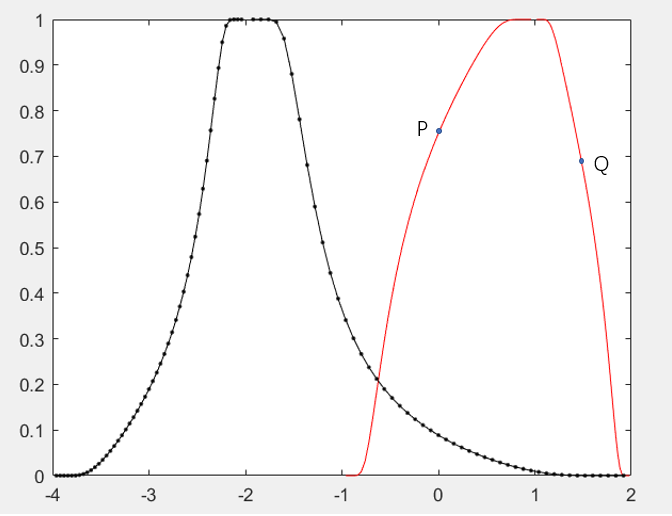}
		\caption{The graph of $(-2)\tilde{b}_{3}$}\label{Fig.4} 
	\end{minipage}
\end{figure}

Figure \ref{Fig.3} (resp. \ref{Fig.4}) shows the fuzzy scalar multiplication with opposite (resp. negative) $\lambda$. When $\lambda > 0$,  the membership value turns to be larger than the original one. Moreover,  the support of the fuzzy number spreads and leads to more uncertainty. For the case $\lambda < 0$,  the graph of resulting fuzzy number is mirror flipped  and shifted to the left horizontally, which is in accordance with our definition in \ref{def}.

\textbf{Example 3.} (Subtraction) 
Set $\tilde{b}_{1},\tilde{b}_{3}$ be as above.  The Gaussian PDMF resulting from fuzzy subtraction is calculated using the operation law $(3)$ in Definition \ref{def}:
	\begin{equation*} 
	\tilde{b}_{3}\ominus\tilde{b}_{1}=\langle(-1,1,2);-0.6745,-0.4399 \rangle\ominus \langle(-1,0,1);0,0\rangle =\langle(-2,1,3);-0.6745,-0.4399 \rangle.
	\end{equation*}
\begin{figure}[htbp]\label{Fig.5} 	
	\centering
	\includegraphics[height=6cm,width=8cm]{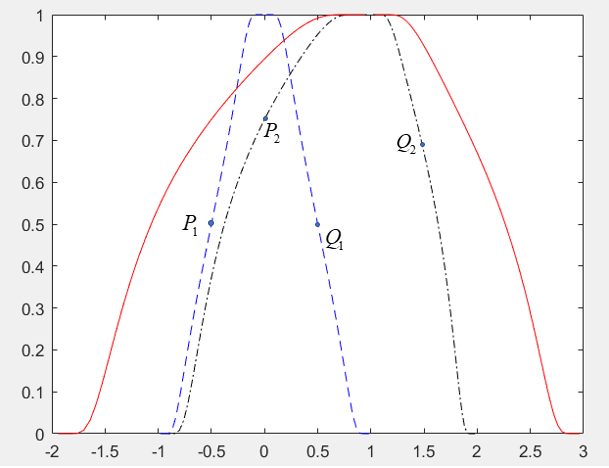}\caption{ Fuzzy subtraction $\tilde{b}_{3}\ominus\tilde{b}_{1}$}
\end{figure}
Figure \ref{Fig.5}  shows the Gaussian PDMF  $\tilde{b}_{3}\ominus\tilde{b}_{1}$.
Note that, the support of the fuzzy number resulting from 
fuzzy subtraction using our proposed method is larger than the support of
MIN and MAX. 

We emphasize that, all fuzzy computations in this examples happen in the same function space, saying,  G-PDMFS.

\section{Final remarks}\label{sec6}

In this paper we presented a new class of fuzzy numbers $X_{h,p}$ in which each fuzzy number is uniquely identified by a membership function $f(x)$ with the form \eqref{fxabc}. More precisely, $f(x)$ is constructed by combining a class of nonlinear mapping $h$ (see Definition \ref{h-} and Definition \ref{h+}) and a class of probability density function $p$ (See Definition \ref{pdf}). Here $h$ can be seen as the subjective perception and $p$ as the objective entity, respectively. The existence of the pair $(h, p)$  is shown for any pre-given information $(a,b,c;P, Q)$ of the fuzzy number.  Especially, the common triangular number  can also be interpreted by a function pair $(h,p)$. 

Next we consider a sample function space $X_{h,p}$ with $h$ being the tangent function and $p$ being the Gaussian-type function with free variable $\mu$. We define the arithmetic operations on $X_{h,p}$ via the free variable $\mu$ which is the expectation of $p(x;\mu)$. Under our definitions,  $X_{h,p}$ has a linear algebra. 

Finally, we provide some numerical examples and graphs of the proposed addition, scalar multiplication and subtraction on the PDMF space $X_{h,p}$.

%\bibliography{/Users/zhengchuang/SynologyDrive/latex/BibTex/ref-fuzzy}
%\bibliography{/Users/chuang/Desktop/SynologyDrive/latex/BibTex/reference-wh}

\end{document}